\documentclass[11pt]{article}
\usepackage{mathrsfs}
\usepackage{amssymb,amsfonts,amsmath,color,latexsym,epsfig,cite,psfrag,eepic,colordvi}
\usepackage{amscd,graphics}
\usepackage{verbatim}
\usepackage{lineno}
\newcommand{\qed}{\hfill $\Box $}
\newcommand{\pf}{\noindent {\bf Proof.} }
\newcommand{\E}{{\mathbb{E}}}
\renewcommand{\P}{{\mathbb{P}}}

\newtheorem{theorem}{Theorem}[section]
\newtheorem{lemma}[theorem]{Lemma}

\newtheorem{definition}[theorem]{Definition}

\newtheorem{observation}[theorem]{Observation}

\begin{document}

\title{Rainbow perfect matchings in 3-partite 3-uniform hypergraphs}

\author{Hongliang Lu\footnote{Partially supported by National Key Research and Development Program of China under grant 2023YFA1010203 and National Natural
Science Foundation of China under grant 12271425}\\
School of Mathematics and Statistics\\
Xi'an Jiaotong University\\
Xi'an, Shaanxi 710049, China\\
\medskip \\
Yan Wang\footnote{Partially supported by National Natural Science Foundation of China under grant 12201400, National Key Research and Development Program of China under grant 2022YFA1006400, and Explore X project of Shanghai Jiao Tong University}\\
School of Mathematical Sciences\\
Shanghai Frontier Science Center of Modern Analysis\\
Shanghai Jiao Tong University\\
Shanghai 200240, China}

\date{}

\maketitle

\date{}

\maketitle

\begin{abstract}
Let  $m,n,r,s$ be nonnegative integers such that $n\ge m=3r+s$ and $1\leq s\leq 3$. Let
\[\delta(n,r,s)=\left\{\begin{array}{ll} n^2-(n-r)^2 &\text{if}\
s=1 , \\[5pt]
 n^2-(n-r+1)(n-r-1) &\text{if}\
s=2,\\[5pt]
n^2 - (n-r)(n-r-1) &\text{if}\ s=3.
\end{array}\right.\]
We show that there exists a constant $n_0 > 0$ such that if $F_1,\ldots, F_n$ are 3-partite 3-graphs with $n\ge n_0$ vertices in each partition class and the minimum vertex
degree of $F_i$ is at least $\delta(n,r,s)+1$ for $i \in [n]$, then $\{F_1,\ldots,F_n\}$ admits a rainbow perfect matching. This generalizes a result of Lo and Markstr\"om on the vertex degree condition for the existence of  perfect matchings in 3-partite  3-graphs.
%In this proof, we use a fractional rainbow matching theory obtained by \textit{Aharoni et al.} to find edge-disjoint fractional perfect matching.

\end{abstract}

\section{Introduction}
Let $k$ be a positive integer. Let $[k]:=\{1,\ldots, k\}$. For a set
$S$, let ${S\choose k}:=\{T\subseteq S: |T|=k\}$. A hypergraph $H$
consists of a vertex set $V(H)$ and an edge set $E(H)$, whose members
are subsets of $V(H)$. A hypergraph $H$ is {\it $k$-uniform} (or a
{\it $k$-graph}) if $E(H)\subseteq {V\choose k}$. Thus, 2-graphs are
precisely graphs without loops.

A {\it matching} in a hypergraph $H$ is a set of pairwise disjoint
edges of $H$, and it is \emph{perfect} if the union of all edges in
the matching is $V(H)$. For a hypergraph $H$, we use $\nu(H)$ to
denote the largest size of a matching in $H$. A 2-graph $G$ is \emph{factor-critical} if  for any $v\in V(G)$, $G-v$ has a perfect matching.

Let $H$ be a $k$-graph and $T\subseteq V(H)$.
The \emph{degree} of $T$ in $H$, denoted by $d_H(T)$, is the number of edges of $H$
containing $T$. Let $l$ be a positive integer.  Let $\delta_l(H):=
\min\{d_H(T): T\in {V(H)\choose l}\}$ be the \emph{minimum
$l$-degree} of $H$. $\delta_1(H)$ is often called the minimum {\it
vertex degree} of $H$, and $\delta_{k-1}(H)$ is known as the minimum
\textit{codegree} of $H$. R\"{o}dl, Ruci\'{n}ski and Szemer\'{e}di
\cite{RRS4} determined the minimum codegree threshold for the
existence of a perfect matching in $k$-graphs. An analogue
for $k$-partite $k$-graphs was obtained by Aharoni, Georgakopoulos
and Spr\"{u}ssel \cite{AGS}.
%Some results are known about the minimum $l$-vertex degree conditions for the appearance of a perfect matching.
H\`{a}n, Person and Schacht \cite{HPS} showed that $\delta_1(H)>(\frac 5 9+o(1)){|H|\choose 2}$
is sufficient for the appearance of a perfect matching in 3-graph $H$, and
they \cite{HPS} further conjectured that the $l$-degree threshold for the existence of a perfect matching in $k$-graphs
is $\delta_l(H)> (\max \left\{\frac 1 2, 1-(1-\frac 1 k)^{k-l}\right\}+o(1)){n\choose {k-l}}$ for $k\geq 3$ and $1\leq l<k$.
K\"{u}hn, Osthus and Treglown \cite{KOT} determined the minimum vertex degree threshold exactly
for the existence of a perfect matching in all large 3-graphs.

Bollob\'{a}s, Daykin  and Erd\H{o}s \cite{BDE} considered the minimum vertex degree for the appearance of matching of size $m$.
%\begin{theorem} (Bollob\'{a}s, Daykin  and Erd\"{o}s \cite{BDE})
They proved that for integers  $k\geq 2$, if $H$ is a $k$-graph of order $n\geq 2k^3m$ and
$\delta_1(H)>{{n-1}\choose {k-1}}-{{n-m}\choose {k-1}},$
then $\nu(H)\ge m$.
%contains a matching of size at least $d$.
%\end{theorem}
For 3-graphs, K\"{u}hn, Osthus and Treglown \cite{KOT} proved a stronger result:
%\begin{theorem} (K\"{u}hn, Osthus and Treglown \cite{KOT})
There exists an $n_0\in \mathbb{N}$ such that if $H$ is a 3-graph of order $n\geq n_0$, $m\le n/3$ and
$\delta_1(H)>{{n-1}\choose 2}-{{n-m}\choose 2},$
then $\nu(H)\ge m$.
%Recently, there have been much effort to extend the minimum degree conditions to rainbow matchings, see for example, \cite{CHWW2023, LWY2022,LWY2023jcta,LWY2023jctb,LWY2024}.
%\end{theorem}

In this paper, we study matchings in $k$-partite $k$-graphs. This is a natural generalization of matchings in
bipartite graphs.  Given a positive integer $k$, a $k$-graph is called $k$-$partite$ if there exists a
partition of the vertex set $V(H)$ into $k$ sets $V_1, V_2, \cdots, V_k$ (called {\it classes}) such that for any $f\in E(H)$, $|f\cap V_i|=1$
for $i\in [k]$.

%Matchings in $k$-partite $k$-graphs also draw much attention recently.
Minimum codegree thresholds for the existence of perfect matchings and near perfect matchings in $k$-partite $k$-graphs have been determined in \cite{LWY2019} and \cite{LWY2018} respectively.
For a
$k$-partite $k$-graph $H$, a set $T\in V(H)$ is said to be {\it
legal} if $|T\cap V_i|\leq 1$ for all $i\in [k]$. Thus, if $T$ is
not legal in $H$ then $d_H(T)=0$.
We define the minimum $l$-degree of $H$ to be $\delta_l(H):= \min\{d_H(T):
T\in {V(H)\choose l} \mbox{ and $T$ is legal}\}$. Lo and
Markstr\"om~\cite{LM} defined the following family of $k$-partite
$k$-graphs. For positive integers $k,l,n$ and nonnegative integers
$d_i$, $i\in [k]$, let $H_{k,l}(n;d_1,d_2,\ldots,d_k)$ denote the
$k$-partite $l$-graph with partition classes $V_1, \ldots, V_k$ such
that for $i\in [k]$, $|V_i|=n$, each $V_i$ has a partition $U_i,W_i$ with
$|W_i|=d_i$, and $E(H_{k,l}(n;d_1,d_2,\ldots,d_k))$ consists of all edges intersecting
$\cup_{i\in [k]}W_i$.
let $H_{k,l}'(n;d_1,d_2,\ldots,d_k)$ denote the
$k$-partite $l$-graph with partition classes $V_1, \ldots, V_k$ such
that for $i\in [k]$, $|V_i|=n$, each $V_i$ has a partition $U_i,W_i$ with
$|W_i|=d_i$, and $E(H_{k,l}'(n;d_1,d_2,\ldots,d_k))$ consists of all edges intersecting
both $\cup_{i\in [k]}W_i$ and $\cup_{i\in [k]} V_i \setminus W_i$.
In particular, when $k=l$, we denote $H_{k,k}(n;d_1,d_2,\ldots,d_k)$ by $H_{k}(n;d_1,d_2,\ldots,d_k)$, and denote $H_{k,k}'(n;d_1,d_2,\ldots,d_k)$ by $H_{k}'(n;d_1,d_2,\ldots,d_k)$.
If $d_i\in \{\lceil m/k\rceil,\lfloor
m/k\rfloor\}$ and $\sum_{i=1}^k d_i=m$, let
$H_k(n,m):=H_k(n;d_1,d_2,\ldots,d_k)$ and $H_k'(n,m):=H_k'(n;d_1,d_2,\ldots,d_k)$.
Define ${\cal H}_k(n;m;d)$ to
be the family of all $k$-partite $k$-graphs $H(n;d_1,\ldots, d_k)$
with $\max\{d_i:i\in [k]\}=d$ and $\sum_{i\in [k]}d_i=m$.
 Define ${\cal H}_k^*(n;m)$ to be the family of all $k$-partite $k$-graphs $H_k(n;m-1)\cup H'$, where $V(H')=V(H_k(n;m-1))$ and
$E(H')$ is an intersecting family.

Lo and Markstr\"om \cite{LM} proved the following.

\begin{theorem} [Lo and Markstr\"om \cite{LM}]\label{thm:LMk}
Let $m,n,k$ be nonnegative integers with $k\ge 2$ and $n\ge k^7m$, and let  $H$ be a $k$-partite $k$-graph with $n$ vertices in each class.
If $\nu(H)= m$ and
$$\delta_1(H)\geq\max\{ \delta_1(F)\ |\ F\in {\cal H}_k(n;m;\lceil m/k\rceil)\},$$
then $H$ is a subgraph of some member of ${\cal H}_k(n;m;\lceil m/k\rceil)\cup {\cal H}_k^*(n;m)$. Moreover,
if $m\ne 1 \pmod k$ then $H$ is a subgraph of a member of ${\cal H}_k(n;m;\lceil m/k\rceil)$.
\end{theorem}

They also computed $\delta_1({\cal H}_3(n;m;\lceil m/k\rceil)$ explicitly and prove the following theorem for $3$-partite $3$-graphs.
 \[d_3(n,m)=\left\{\begin{array}{ll} n^2-(n-\lfloor m/3\rfloor)(n-\lfloor (m+1)/3\rfloor) &\text{if}\
m\ne 1 \pmod 3, \\[5pt]
 n^2-(n-(m-1)/3)^2+1 &\text{if } m=1 \pmod 3.\\[5pt]
\end{array}\right.\]

\begin{theorem} [Lo and Markstr\"om \cite{LM}]\label{thm:LM3}
If $H$ is a $3$-partite $3$-graph with each class of size $n\ge
4\times3^6m$ and $\delta_1(H)> d_3(n,m)$ then $\nu(H)\ge m+1$.
\end{theorem}

For perfect matchings, this bound can be further extended to graphs with sufficiently large $n$: There
is an integer $n_0$ such that if $H$ is a $3$-partite $3$-graph with
$n\ge n_0$ vertices in each class and $\delta_1(H)> d_3(n,n-1)$ then
$H$ has a perfect matching.

Lo and Markstr\"om \cite{LM} asked whether $\nu(H)\ge m+1$ for every
$3$-partite $3$-graph $H$ with each class of size $n>m$ and
$\delta_1(H)>d_3(n,m)$.
Lu and Zhang \cite{LZ17} answered Lo and Markstr\"om's problem by showing that this is true
provided $n\ge n_0$ for some constant $n_0$.

%In this paper, we show that this is true
%provided $n\ge n_0$ for some constant $n_0$. This is done by
%refining and combining some techniques from K\"{u}hn, Osthus and
%Treglown \cite{KOT} and Lo and Markstr\"om \cite{LM}, which in turn
%was inspired by  H\`{a}n, Person and Schacht \cite{HPS}.

\begin{theorem}[Lu and Zhang \cite{LZ17}]\label{LM-MAIN}
Let  $m,n,r,s$ be nonnegative integers such that  $m=3r+s$ and
$1\leq s\leq 3$. Let
\[\delta(n,r,s)=\left\{\begin{array}{ll} n^2-(n-r)^2 &\text{if}\
s=1 , \\[5pt]
 n^2-(n-r+1)(n-r-1) &\text{if}\
s=2,\\[5pt]
n^2 - (n-r)(n-r-1) &\text{if}\ s=3.
\end{array}\right.\]
There exists an $n_0\in \mathbb{N}$ such that if $H$ is a
$3$-partite $3$-graph with $n\ge n_0$ vertices in each class and
$\delta_1(H)>\delta(n,r,s)$ then $\nu(H)\ge m$.
\end{theorem}

Let $F_1,\ldots, F_t$ be $t$ hypergraphs and let
$\mathcal{F} = \{F_1,\ldots, F_t\}$  denote a family
of hypergraphs.
%when there is no confusion;
A set of pairwise disjoint edges, one from
each $F_i$, is called a \emph{rainbow matching} for $\mathcal{F}$. In
this case, we also say that ${\cal F}$ {\it admits} a rainbow
matching.
There has been much effort to generalize minimum degree conditions for matchings in hypergraphs to rainbow matchings in a family of hypergraphs, see \cite{AHJ,CHWW2023,LWY2022,LWY2023jcta,LWY2023jctb,LWY2024}.

In this paper, we first show the following rainbow version of Theorem~\ref{thm:LMk} when $k = 3$.

\begin{theorem}\label{small-matchings}
Let  $n,r,s$ be nonnegative integers such that  $m=3r+s$ and
$1\leq s\leq 3$. Let
\[\delta(n,r,s)=\left\{\begin{array}{ll} n^2-(n-r)^2 &\text{if}\
s=1 , \\[5pt]
 n^2-(n-r+1)(n-r-1) &\text{if}\
s=2,\\[5pt]
n^2 - (n-r)(n-r-1) &\text{if}\ s=3.
\end{array}\right.\]
Suppose $n \ge 120 m$.
If $\{F_1,\ldots, F_{m}\}$ is a family of
$3$-partite 3-graphs on the same vertex set  and with $n$ vertices in each partition class  and
$\delta_{1}(F_i) >\delta(n,r,s)$ for $i\in [m]$.
Then $\{F_1,\ldots, F_{m}\}$ admits a rainbow matching of size $m$.
\end{theorem}

Moreover, we also determine the exact minimum vertex degree condition for the existence of rainbow perfect matchings in a family of $3$-partite $3$-graphs.

\begin{theorem}\label{main}
Let  $n,r,s$ be nonnegative integers such that  $n=3r+s$ and
$1\leq s\leq 3$. Let
\[\delta(n,r,s)=\left\{\begin{array}{ll} n^2-(n-r)^2 &\text{if}\
s=1 , \\[5pt]
 n^2-(n-r+1)(n-r-1) &\text{if}\
s=2,\\[5pt]
n^2 - (n-r)(n-r-1) &\text{if}\ s=3.
\end{array}\right.\]
There exists an $n_0\in \mathbb{N}$ such that the following result holds. If $\{F_1,\ldots, F_{n}\}$ is a family of
$3$-partite $3$-graphs on the same vertex set  and with $n\geq n_0$ vertices in each partition class  and
$\delta_{1}(F_i) >\delta(n,r,s)$ for $i\in [n]$.
Then $\{F_1,\ldots, F_{n}\}$ admits a rainbow perfect matching.
\end{theorem}

It is easy to see that we can derive the perfect matching result of Theorem~\ref{LM-MAIN} from
Theorem~\ref{main} by setting $F_1=\ldots =F_{n}=H$.  Moreover, if $F_i=H(n,k)$ for $i\in [n]$ then $\{F_1, \ldots,
F_{n}\}$ admits no rainbow perfect matching.  So the vertex degree condition in Theorem~\ref{main} is best possible.

In the rest of this section, we give some notations that will be used in this paper.
By $x\ll y$, we mean that $x$ is sufficiently small compared with $y$ which need to satisfy finitely many inequalities in the proof.
Let  $Q=\{v_1,\ldots,v_{n}\}$ and $\mathcal{F}=\{F_1,\ldots,F_{n}\}$ be a family of $3$-partite $3$-graphs on the same vertex set. Let  $H_{4}(\mathcal{F})$ be a balanced $4$-partite 4-graph with vertex partition $Q \cup V(F_1)$ and edge set $E(H_{4}(\mathcal{F}))=\cup_{i=1}^{n}E_i$, where $E_i=\{e\cup \{v_i\}\ |\ e\in E(F_i)\}$ for $1\leq i\leq n$.
In particular, if  $F_1,\ldots,F_{n}$ are $n$ copies of $H_3(n,n-1)$, then we write $H_{1,3}(n):=H_{4}(\mathcal{F})$.
Furthermore, if  $F_1,\ldots,F_{n}$ are $n$ copies of $H_3'(n,n-1)$, then we write $H_{1,3}'(n):=H_{4}(\mathcal{F})$.
It is easy to see the following.
\begin{observation}\label{Obser1}
${\cal F}$ has a rainbow matching if and only if $H_{4}({\cal F})$ has a perfect matching.
 \end{observation}

We also need the following definition. Let \(\varepsilon > 0\) be a real number. We say that \(H_2\) is $\varepsilon$-close of $H_1$ if the following condition holds:
\[
c(H_1, H_2) \leq \varepsilon |V(H_1)|^k,
\]
where $c(H_1, H_2)$ represents the smallest possible size of the set difference $|E(H_1) \setminus E(H_2')|$, considering all graphs $H_2'$ that are isomorphic to $H_2$ and share the same vertex set as $H_1$.

The organization of this paper is as follows.
In Section 2, we prove Theorem \ref{small-matchings} for the existence of rainbow matchings of size $m$ when $n \ge 120 m$.
By Observation \ref{Obser1}, it suffices to show that $H_{4}({\cal F})$ has a perfect matching.
We divide the proof by considering whether $H_{4}({\cal F})$ is close to the extremal graph $H_{1,3}'(n)$ or not.
We show the case when it is close to the extremal graph in Section 3.
In Section 4, we prove an absorbing lemma which is needed in the proof of the case when $H_{4}({\cal F})$ is not close to the extremal graph.
We show the existence of an almost perfect matching in Section 5.
Finally, we prove Theorem \ref{main} in Section 6.

\section{Small rainbow matchings}

In this section, we prove Theorem \ref{small-matchings}.
First, we need to show the existence of a rainbow matching of size two in two $3$-partite $3$-graphs when minimum vertex degree is at least two.

\begin{lemma} \label{small-matchings-size-two}
Let $n_1,n_2,n_3$ be three integers such that $\max\{n_1,n_2,n_3\}\leq 3\min\{n_1,n_2,n_3\}/2$.
Let $F_1,F_2$ be two $3$-partite $3$-graphs on the same vertex partition classes $V_1,V_2,V_3$, where $n_i=|V_i|$ for $i\in [3]$. If  for $i\in [2]$, $\delta_1(F_i)>1$, then $\{F_1,F_2\}$ admits a rainbow matching.
\end{lemma}
\pf Let $e\in E(F_1)$. If $E(F_2-V(e))\neq \emptyset$, let $e'\in E(F_2-V(e))$, then $\{e',e\}$ is a rainbow matching of $\{F_1,F_2\}$. So we may assume every edge of $E(F_2)$ intersecting $e$. Since $\delta_1(F_2)\geq 2$, then there exists  $x\in e\cap (V_2\cup V_1)$ such that $d_{F_2}(x)\geq n_3$. Without loss generality, suppose that $x\in V_1$. Since $\delta_1(F_1)\geq 2$, $E(F_1-x)\neq\emptyset $. We choose $g\in E(F_1-x)$. If $F_2-V(g)$ has an edge $g'$, then $\{g,g'\}$ is a desired rainbow matching. So we may assume that every edge of $F_2$ intersects $g$. So there exists $y\in g\backslash V_1$ such that $d_{F_2}(y)\geq \frac{2(n_1 - 1)+n_3}{2} = n_1+n_3/2-1$.

 If $N_{F_2}(y)$ has a matching of size three, say $M=\{f_1,f_2,f_3\}$. 
Since  $E(F_1 - y) \ne \emptyset$, then there exists $e'\in E(F_1-y)$.
Note that there exists $f\in M$ such that $f\cap e'=\emptyset$. Then $\{f\cup \{y\},e'\}$ is a rainbow matching of $\{F_1,F_2\}$.

By K\"{o}nig Theorem, $N_{F_2}(y)$ has a vertex cover of size two, say $\{u,v\}$.
Since $d_{F_2}(y)\geq n_1+n_3/2-1$, we have that
 $N_{F_2}(y)$ has two vertex disjoint paths $P_1=v_1vv_2$ and $P_2=u_1u$. Let $w\in V_1$. Since $\delta(F_1)\geq 2$, there exists an edge $h\in E(F_1)$ such that $w\in h$ and $|h\cap \{u,v\}|\leq 1$. So one can see that there exists $f\in\{vv_1,vv_2,uu_1\}$ such that $f\cap h=\emptyset$. Thus $\{h,f\cup\{y\}\}$ is the desired rainbow matching.  This completes the proof. \qed

Next, we prove the existence of a rainbow perfect matching in three $3$-partite $3$-graphs given vertex degree conditions.

\begin{lemma}\label{small-matchings-size-3}
Let $n_1,n_2,n_3$ be three integers such that $n_3=\max\{n_1,n_2,n_3\} \leq \frac{3}{2}\min \{n_1,n_2,n_3\}$.
Let $F_1,F_2,F_3$ be three $3$-partite $3$-graphs with the same vertex classes $V_1,V_2,V_3$, where $n_i=|V_i|$ for $i\in [3]$. For $j\in [3]$, if the follow two conditions hold, then $\{F_1,F_2,F_3\}$ admits a rainbow matching.
\begin{itemize}
  \item $d_{F_j}(v)>n_3$ for $v\in V_1\cup V_2$;
  \item  $d_{F_j}(u)>n_2$ for $u \in V_3$.
\end{itemize}
\end{lemma}
\pf By Lemma \ref{small-matchings-size-two}, $\{F_1,F_2\}$ admits a rainbow matching, say $M=\{e_1,e_2\}$. If $E(F_3\backslash V(M))\neq \emptyset$, then there exists $e_3\in E(F_3\backslash V(M))$ such that $\{e_1,e_2,e_3\}$ is a rainbow matching of $\{F_1,F_2,F_3\}$.   Since $d_{F_3}(v)>n_2$ for all $v\in V_3$, then there exists $x\in V(M)\cap (V_1\cup V_2)$ such that $d_{F_3}(x)\geq n_3^2/8$. Without loss generality, suppose that $x\in V_1$.

For $i\in \{1,2\}$, put $F_i'=F_i-x$. Let $f_1\in F_1'$. If there exists $f_2\in F_2'$ such that $\{f_1,f_2\}$  is  matching, then we may choose an edge $f_3$ in $F_3$ containing $x$ such that $\{f_1,f_2,f_3\}$ is a rainbow matching of $\{F_1,F_2,F_3\}$. So we may assume that every edge of $F_2'$ intersects $f_1$. Thus there exists $y\in V_2\cup V_3$ such that $d_{F_2'}(y)>n_3^2/8$. By the degree condition, one can see that $E(F_1-x-y)\neq \emptyset$. Let $e_1'\in F_1-x-y$. Since $d_{F_2'}(y)>n_3^2/8$, we may pick an edge $e_2'$ in $F_2'-V(e_1')$ containing $y$. Similarly, we may pick an edge $e_3'$ from $F_3-V(\{e_1',e_2'\})$ containing  $x$. Then $\{e_1',e_2',e_3'\}$ is a desired rainbow matching. This completes the proof. \qed

Now we are ready to show Theorem \ref{small-matchings} by induction.

\noindent \textbf{Proof of Theorem \ref{small-matchings}.}
By induction on $n+m$.
By Lemmas \ref{small-matchings-size-two} and \ref{small-matchings-size-3}, we may assume that $m\geq 4$.
Suppose that the result holds for smaller $n+m$ such that $n\geq 120m$.
First, consider that there exists $i\in [m]$, such that $F_i$ contains a matching $N_1$ of size $3m-2$.
Without loss of generality, we may assume $i=m$.
By inductive hypothesis,
$F_1,\ldots, F_{m-1}$ containing a rainbow matching, say $M_2$. Note that $|V(M_2)|= 3m-3$. Thus $N_1$ contains an edge $e$ such that $e\cap V(M_2)=\emptyset$. Then $M_2\cup\{e\}$ is a desired rainbow matching of $\{F_1,\ldots, F_{m}\}$. So we may assume $\nu(F_i)\leq 3m-3$ for all $i\in [m]$.

\medskip
\textbf{Claim 1.} If  there exist $\{i_1,i_2,i_3\}\in {[m]\choose 3}$, and $v_1\in V_1, v_2\in V_2, v_3\in V_3$ such that $d_{F_{i_j}}(v_j)>2mn$ for $j\in [3]$, then $\{F_1,\ldots,F_m\}$ admits a rainbow matching.
\medskip

Suppose the result does not hold. Without loss generality, we may assume  that there exist  $v_1\in V_1, v_2\in V_2, v_3\in V_3$ such that $d_{F_{j}}(v_j)>2mn$. Then for any $i\in [3]$,
\[\delta_1(F_i - \{v_1,v_2,v_3\})>\left\{\begin{array}{ll} (n-1)^2-(n-r)^2 &\text{if}\
s=1 , \\[5pt]
 (n-1)^2-(n-r+1)(n-r-1) &\text{if}\
s=2,\\[5pt]
(n-1)^2 - (n-r)(n-r-1) &\text{if}\ s=3.
\end{array}\right.\]
Thus by inductive hypothesis, $\{F_i\setminus \{v_1,v_2,v_3\}\ |\ i\in [m] \backslash [3]\}$ contains a rainbow matching $M_1$. Since  $d_{F_{1}}(v_1)>2mn$ and the edges in $F_1$ containing $v_1$ and intersecting $V(M_1)\cup \{v_2,v_3\}$ is at most $2nm-m^2$, we can find an edge $e_1$ in $F_1-V(M_1)-\{v_1,v_2,v_3\}$ containing $v_1$.
Similarly, we can greedily find $e_2\in F_{2} - V(M_1) - \{v_1,v_3\} - V(e_1), e_3\in F_{3} - V(M_1) - \{v_1,v_2\} - V(\{e_1, e_2\})$ such that $M_1\cup \{e_1,e_2,e_3\}$ is a desired rainbow matching.
This completes the proof of Claim 1.

For $i\in [m]$, denote $A_i:=\{x\in V(F_i)\ |\ d_{F_i}(x)>2mn\}$.

 \medskip
\textbf{Claim 2.} The following statements are true:
\\
(i) $A_i\cap (V_l\cup V_j)\neq \emptyset$ for any $\{l,j\}\in {[3]\choose 2}$;
\\
(ii) If there exists $i\in [m]$ such that   $A_i\setminus V_j\neq \emptyset$ for all $j\in [3]$, then $\{F_1,\ldots,F_m\}$ admits a rainbow matching.
\medskip

 It suffices to show that $A_i\cap (V_1\cup V_2)\neq \emptyset$. Let $M$ be a maximal matching of $F_i$. Then $|M|\leq 3m$.
By inductive hypothesis, $\{F_1,\ldots,F_m\} \setminus \{F_i\}$ admits a rainbow matching of size $m-1$.
   Since  $\delta(F_i)\geq 2nr-r^2$ and $n\geq 100m$, then by the maximality of $m$, the number of edges intersecting $V_3-V(M)$ and $V(M)-V_3$ is at least $(2nr-r^2+1)(n-3m)$. Thus there exists $x\in V(M)\cap (V_1\cup V_2)$ such that
 \[
 (2nr-r^2+1)(n-3m)/(2 \cdot 3m)\geq 14nr-7r^2+14>13nr>6n(r+1)\geq 2nm.
 \]
 This completes the proof of Claim 2 (i).

Suppose that there exists $i\in [m]$ such that  $A_i\cap V_j\neq \emptyset $ for all $j\in [3]$. Let $v_j\in A_i\cap V_j$ for $j\in [3]$.  Let $\{i',i''\}\in {[m]\setminus \{i\}\choose 2}$.
Recall by Claim 2 (i) that  there exist $S,T\in {[3]\choose 2}$ such that $A_{i'}\cap V_j\neq \emptyset$ for all $j\in S$ and $A_{i'}\cap V_j\neq \emptyset$ for all $j\in T$. Without loss generality, we may suppose that $S=\{1,2\}$ and $1\in T$.  Let $u_2\in A_{i'}\cap V_2$ and $w_1\in A_{i''}\cap V_1$.
Then we have $d_{F_i}(v_3) > 2mn$,  $d_{F_i'}(u_2) > 2mn$ and  $d_{F_i''}(w_1) > 2mn$.
%$d_H(x)>2mn$ for any $x\in \{w_1,u_2,v_3\}$.
By Claim 1, $\{F_1,\ldots, F_m\}$ admits a rainbow matching. This completes the proof of Claim 2 (ii).

\medskip

By Claim 2, we may assume that $A_1\cap V_j\neq \emptyset$ for $j\in [2]$.

\medskip
\textbf{Claim 3.~}   If there exists $i\in [m]\backslash \{1\}$ such that  $A_i\cap V_3\neq \emptyset$, then $\{F_1,\ldots,F_m\}$ admits a rainbow matching.
\medskip

Without loss generality, we may assume that $A_2\cap V_2\neq \emptyset$ and $A_2\cap V_3\neq \emptyset$. Let $u_i\in A_1\cap V_i$ for $i\in [2]$ and $v_j\in A_2\cap V_j$ for $j\in \{2,3\}$. Recall that $A_3\cap (V_1\cup V_3)\neq \emptyset$ by Claim 2 (i). By symmetry, we may assume that $A_3\cap V_1\neq \emptyset.$ Let $w_1\in A_3\cap V_1$.
Then we have  $d_{F_3}(w_1) > 2mn$,  $d_{F_2}(v_3) > 2mn$ and  $d_{F_1}(u_2) > 2mn$.
%$d_H(x)>2mn$ for any $x\in \{w_1,u_2,v_3\}$.
Thus by Claim 1, $\{F_1,\ldots,F_m\}$ admits a rainbow matching. This completes the proof of Claim 3.

%By symmetry and Claim 2, we may assume that there exists $i\in [m]\setminus \{1\}$, say $i=2$ such that $A_i\cap V_3=\emptyset$ for all $j\in [3]$. By Claim 2, we have $A_3\cap (V_1\cup V_2)\neq \emptyset$. Without loss generality, suppose $A_3\cap V_1\neq \emptyset$. Now we may choose $v_1\in A_3\cap V_1$, $v_2\in A_1\cap V_2$ and $v_3\in A_2\cap V_3$.
%Then we have $d_H(x)>3mn$ for any $x\in \{v_1,v_2,v_3\}$. Thus by Claim 1, $\{F_1,\ldots,F_m\}$ admits a rainbow matching. This completes the proof of Claim 3.

 %By Claim 2, we may assume there exists $i,j$, say $i=1,j=2$ such that $A_1\cap V_1=\emptyset$ and $A_2\cap V_2=\emptyset.$
%We choose $v_2\in A_1\cap V_2$, $v_3\in V_3\cap A_1$, $u_1\in A_2\cap V_2$, $u_3\in V_3\cap A_2$, where $v_2,v_3,u_1,u_3$ are not necessarily different. Consider $A_3\neq \emptyset$. Without loss generality, we may assume that $A_3\cap V_1\emptyset$. Let $w_1\in A_3\cap V_1$.

 %By Claim 2, $|A_i\cap (V_2\cup V_3)|= |A_i\cap V_2|\geq m/2$ and $|A_i\cap (V_1\cup V_3)|= |A_i\cap V_1|\geq m/2$. So we have the following Claim.

\medskip

By Claim 3, we may assume that $A_i\cap V_3=\emptyset$ for all $i \in [m]$.

\medskip
\textbf{Claim 4.~}  If $m\neq 5$, then $|A_i\cap V_1|\geq m/2$ and $|A_i\cap V_2|\geq m/2$ for all $i\in [m]$; else $|A_i\cap (V_1\cup V_2)|\geq 4$ for {all $i\in [5]$.
\medskip

By contradiction, suppose that $|A_i \cap V_1|< m/2$. Let $M$ be a maximum matching of $F_i$. Note every edge in $F_i$ intersect $V(M)$. Moreover, for every $x\in V(M)\cap V_3$, we have $d_{F_i}(x)\leq 2nm$. Thus
we have
\[
\sum_{x\in V_2-V(M)}d_{F_i}(x)\leq |A_i\cap V_1|n^2+(6m-|A_i|)2nm.
\]
So if $m$ is even, then
\begin{align}\label{m-even}
\sum_{x\in V_2-V(M)}d_{F_i}(x)\leq (m/2-1)n^2+(5m/2+1)2nm.
\end{align}
and else
\begin{align}\label{m-odd}
\sum_{x\in V_2-V(M)}d_{F_i}(x)\leq ((m-1)/2)n^2+(5m/2+1/2)2nm.
\end{align}
On the other hand,
\begin{align}\label{upp-bound}
\sum_{x\in V_2-V(M)}d_{F_i}(x)\geq (n-3m)(\delta(n,r,s)+1)=\left\{\begin{array}{ll} (2nr-r^2+1)(n-3m) &\text{if}\
s=1 , \\[5pt]
(2nr-r^2+2)(n-3m) &\text{if}\
s=2,\\[5pt]
(2nr-r^2+n-r) &\text{if}\ s=3.
\end{array}\right.
\end{align}

Now we divide into four cases according to the values of $s$ and $m$.

\medskip
\textbf{Case 1.~}$s=1$.
\medskip

 By (\ref{m-odd}) and (\ref{upp-bound}),  we have
\begin{align*}
(m-1)n^2/2+(5m+1)nm&>(2nr-r^2+1)(n-3m)\\
&=(2n(m-1)/3-(m-1)^2/9+1)(n-3m)\\
&>(2n(m-1)/3-(m-1)^2/9)(n-3m)\\
&>2n^2(m-1)/3-2nm(m-1)-m(m-1)n/9.
\end{align*}
So we have
\begin{align*}
(m-1)n/2+(5m+1)m>2n(m-1)/3-2m(m-1)-m(m-1)/9,
\end{align*}
which implies that
\begin{align*}
10m(m-1)<n(m-1)/6<(5m+1)m+2m(m-1)+m(m-1)/9=m(7m-1)+m(m-1)/9
\end{align*}
since $n\geq 60m$.
%So we may infer that
%\[
%26m/9<80/9,
%\]
But this is a contradiction since $m\geq 4$.

\medskip
\textbf{Case 2.~}$s=3$ (so $m \ge 6$).
\medskip

By (\ref{m-odd}) and (\ref{upp-bound}),  we have
\begin{align*}
(m-1)n^2/2+(5m+1)nm&>(2nr-r^2+n-r)(n-3m)\\
&=(2n(m-3)/3-(m-3)m/9+n)(n-3m)\\
&=(2n(m-1)/3-n/3-(m-3)m/9)(n-3m)\\
&>2n^2(m-1)/3-2nm(m-1)-n^2/3+nm-(m-3)nm/9.
\end{align*}
So we have
\begin{align*}
(5m+1)m>n(m-1)/6-2m(m-1)-n/3+m-(m-3)m/9,
\end{align*}
which implies that
\begin{align*}
20m^2-60m\leq n(m-3)/6<(5m+1)m+2m(m-1)-m+(m-3)m/9=7m^2+m(m-3)/9-2m
\end{align*}
since $n\geq 120m$.
%Thus we may infer that
%\[
%116m\leq 519,
%\]
But this is a contradiction since $m\geq 6$.

\medskip
\textbf{Case 3.~}$m=5$.
\medskip

Suppose that the claim does not hold. Then there exists $i\in [5]$ such that $|A_i\cap (V_1\cup V_2)|\leq 3$. %and $|V(M)|\leq 12$.
 %Without loss generality, we assume that $i=5$.
Let $M'$ be a maximum matching of $F_i$. Recall that $|M'|\leq 12$. Then we have
\begin{align*}
\sum_{x\in (V_1\cup V_2)-V(M')}d_{F_i}(x)&\leq 3n^2+(|V(M')\cap (V_1\cup V_2)|-3)2nm+2|V_3\cap V(M')|2nm\\
&\leq 3n^2+93nm.
\end{align*}
So we have
\begin{align}\label{m5-odd}
\sum_{x\in (V_1\cup V_2)-V(M')}d_{F_i}(x)\leq 3n^2+93nm.
\end{align}
On the other hand,
\begin{align}\label{upp-m5}
\sum_{x\in (V_2\cup V_1)-V(M')}d_{F_i}(x)\geq 2(n-3m)(\delta(n,r,s)+1)\geq
2(2n+1)(n-3m)>4n^2+n-12nm.
\end{align}
Combining (\ref{m5-odd}) and (\ref{upp-m5}), we have $n+1<105m$,
a contradiction since $n\geq 120m.$

\medskip
\textbf{Case 4.~}$m\geq 8$ and $s=2$.
\medskip

By (\ref{m-odd}) and (\ref{upp-bound}),  we have
\begin{align*}
&(m-1)n^2/2+(5m+1)nm\\
&>(2nr-r^2+2)(n-3m)\\
&=(2n(m-1)/3-2n/3-(m-2)^2/9+1)(n-3m)\\
&>2n^2(m-1)/3-2nm(m-1)-2n^2/3-n(m-2)^2/9+2nm.
\end{align*}
So we have
\begin{align*}
n(m-1)/6-2n/3<2m(m-2)+(m-2)^2/9+(5m+1)m.
\end{align*}
Thus, since $n\geq 120m$, we have
\begin{align*}
20m^2-100m&\leq n(m-1)/6-2n/3\\
&<7m^2-3m+(m-2)^2/9.
\end{align*}
%i.e.,
%\begin{align*}
%13m-97<(m-2)^2/9m<m/9,
%\end{align*}
But this is a contradiction since $m\geq 8$. This completes the proof of Claim 4.
%Since $M$ is a maximum matching, every edge of $H$ intersects $V(M)$. Suppose that $|A_i\cap V_2|<m/2$.
%Then we have
%\[
%(2rn-r^2)(n-m)\leq \sum_{v\in V_1-V(M)}d_{H}(v)<\sum_{v\in (V_2\cup V_3)\cap V(M)}\leq mn^2/2+3mn\cdot 3m/2=mn(n+9m)/2.
%\]
%Note that
%$(2rn-r^2)(n-m)\geq (2mn/3-m^2/9)(n-m)=\frac{1}{9}m(m^2+6n^2-5nm)\geq \frac{2}{3}mn(n-m)$.
%Thus we may infer that
%\[
%\frac{2}{3}mn(n-m)<nm(n+9m)/2,
%\]
%which implies that $n<31m$, a contradiction.

\medskip

Now consider the case when $m=5$.
By Claim 4, $|A_i\cap (V_1\cup V_2)|\geq 4.$
So we can find four different vertices $\{v_1,v_2,v_3,v_4\}$ such that $v_i\in A_i\cap (V_1\cup V_2)$.
Since $\delta_1(F_5)\geq 2n+1$, then $E(F_5)\backslash \{v_1,v_2,v_3,v_4\}\neq \emptyset$. Let $e_5\in E(F_5)\backslash \{v_1,v_2,v_3,v_4\}$. Note that the number of edges in $F_4$ containing $v_4$ and intersecting  $e_5\cup \{v_1,v_2,v_3\}$ is at most $5n$. Since $d_{F_4}(v_4)>2mn=10n$, then we can choose an edge $e_4\in E(F_4-V(e_5)-\{v_1,v_2,v_3\})$.
With similar discussion, we can choose $e_1\in E(F_1),e_2\in E(F_2),e_3\in E(F_3)$ such that $\{e_i\ |\ i\in [5]\}$ is a rainbow matching of  $\{F_i\ |\ i\in [5]\}$.

Finally consider $m\geq  4$ and $m\neq 5$. By Claim 4 and Hall's Theorem, there exist $m$ different vertices, say $v_1,\ldots,v_{m}$ such that $v_i\in A_i \cap (V_1 \cup V_2)$. Write $S_i=\{v_1,\ldots, v_{i}\}$.
One can see that $d_{F_1}(v_1)>2nm$, $d_{F_2-S_1}(v_2)>n(2m-1)$,\ldots, $d_{F_m-S_{m-1}}(v_m)>nm$.
Now we greedily find a rainbow matching. Since  $d_{F_m-S_{m-1}}(v_m)>nm$, we can find an edge $e_1\in F_m-S_{m-1}$ and $v_m\in e_1$. Suppose we have found edges $e_1\in F_m-S_{m-1}$, $e_2\in F_{m-1}-(S_{m-1}\cup e_1)$, \ldots, $e_i\in F_{m-i+1}-(S_{m-i}\cup \bigcup_{j=1}^{i-1} V(e_j))$ and $v_{m-j+1}\in e_j$ for $j\in [i]$. Since $d_{F_{m-i}-S_{m-i-1}}(v_{m-i})>(m+i)n$, then there exists an edge $e_{i+1}$ containing $v_{m-i}$ and avoiding
$S_{m-i-1}\cup (\bigcup_{j=1}^i e_j)$. Continuing this process, we obtain a rainbow matching of $\{F_1,\ldots,F_{m}\}$. \qed

\section{Extremal case}

In this section, we prove the existence of perfect matchings in $H$ when $H$ is close to the extremal graph $H_{1,3}'(n)$.

First, we show when every vertex of $H$ has almost the same neighborhood in $H$ as in the extremal graph $H_{1,3}'(n)$, then $H$ has a perfect matching. To complete the proof,  we need to provide a more precise definition for ``close".
A vertex $v$ of ${\cal F}$ is said to be {\it $\alpha$-good with respect
to} ${\cal H}$ if
$| N_{{\cal H}}(v) \setminus  N_{{\cal F}}(v) | < \alpha |V(\mathcal{F})|^{k}$.
Otherwise, $v$ is said to be {\it
  $\alpha$-bad} with respect to ${\cal H}$.

\begin{lemma}\label{good-PM}
Let  $n$ be nonnegative integers such that $1/n\ll\alpha\ll 1$, $n\equiv 0\pmod 3$.
Let $H$ be a balanced
$4$-partite $4$-graphs with $n$ vertices in each partition class and with partition  $V_1\cup V_2\cup V_3\cup V_4$.
If every vertex in $H$ is $\alpha$-good with respect to $H_{1,3}'(n)$, then $H$ has a perfect matching.
\end{lemma}

\pf For $i\in [4]\backslash\{1\}$, let $U_i\cup W_i$ be a partition of $V_i$ that corresponds to the definition of $H_3(n,n-1)$ and let $W_1=\emptyset, U_1:=V_1$.
Let $M$ be a maximum matching of $H$ such that for every $e\in M$, $|e\cap (W_2\cup W_3\cup W_4)|=1$. Firstly, we show that $|M|\geq n/2$. Otherwise, suppose that $|M|\leq n/2$.
Then $|(W_2\cup W_3\cup W_4)\backslash V(M)|\geq n/2$ and $|(U_2\cup U_3\cup U_4)\backslash V(M)|\geq n$. Without loss generality, suppose that $|W_2\backslash V(M)|\geq |W_3\backslash V(M)|\geq |W_4\backslash V(M)|$. Then $|W_2\backslash V(M)|\geq n/6$, $|U_3\backslash V(M)|\geq n/6$ and $|U_4\backslash V(M)|\geq n/3$.
Let $x\in V_1$. Then
\[
|N_{H_{1,3}'(n)}(x)\backslash N_H(x)|\geq |W_2\backslash V(M)||U_3\backslash V(M)||U_4\backslash V(M)|\geq n^3/108>\alpha n^3,
\]
contradicting to the assumption that $x$ in $H$ is $\alpha$-good with respect to $H_{1,3}(n)$.

Next we show that $M$ is a perfect matching of $H$. By contradiction, suppose that $M$ is not a perfect matching. Without loss generality, suppose
 $x_i\in V_i\backslash (V(M)\cup W_i)$ for $i\in [3]$ and $x_4\in W_4\backslash V(M)$.  Given three edges $E_1,E_2, E_3$ of distinct matching edges
from $M$, we say that $E_1E_2E_3$ is good for $x_1,x_2,x_3,x_4$ if there are all four possible   vertex-disjoint edges $E$ in $H$ which take the
following form: $E$ has type $UUUW$ and contains one vertex from $\{x_1, x_2, x_3,x_4\}$, one vertex from $E_i$ for $1\leq i\leq 3$. Note that if $E_1E_2E_3$ is good for $\{x_1, x_2, x_3,x_4\}$,  then $H$ has a matching of size 4 which consists of
edges of type $UUUW$ and contains precisely the vertices in $E_1,E_2,E_3$ and $\{x_1, x_2, x_3,x_4\}$. So if such a tuple
$E_1E_2E_3$ exists, we obtain a matching in $H$ that is larger than $M$, yielding a contradiction.

So we may assume there does not exists such  $E_1E_2E_3$. Then there exists $x\in \{x_1, x_2, x_3,x_4\}$ such that
$$|N_{H_{1,3}'(n)}(x)\setminus N_H(x)| \geq \frac{1}{4}{|M|\choose 3} >\alpha n^3,$$
contradicting to the assumption that $x$ is $\alpha$-good. This completes the proof. \qed

Now we are ready to show the following main theorem of this section.
We match the bad vertices first and then apply Lemma \ref{good-PM} to find a perfect matching on the remaining good vertices.

\begin{theorem}\label{close-PM}
Let  $n,r,s$ be nonnegative integers such that $0<1/n\ll\varepsilon\ll 1$, $n=3r+s$ and
$1\leq s\leq 3$. Let
\[\delta(n,r,s)=\left\{\begin{array}{ll} n^2-(n-r)^2 &\text{if}\
s=1 , \\[5pt]
 n^2-(n-r+1)(n-r-1) &\text{if}\
s=2,\\[5pt]
n^2 - (n-r)(n-r-1) &\text{if}\ s=3.
\end{array}\right.\]
Let $H$ be a
$4$-partite $4$-graph with $n$ vertices in each partition class and with partition  $V_1\cup V_2\cup V_3\cup V_4$.
If $H$ is $\varepsilon$-close to $H_{1,3}'(n)$  and $d_{H}(\{x,y\}) >\delta(n,r,s)$ for  $x\in V_1$ and $y\in V_2\cup V_3\cup V_4$, then $H$ has a perfect matching.
\end{theorem}

\pf  Suppose that for $2\leq i\leq 4$, let $U_i$  and $W_i$ denote the vertex
classes of $V_i$ of size $n/3$ and $2n/3$ respectively as in the definition of $H_{1,3}'(n)$.
Let $W_1=\emptyset$ and $U_1:=V_1$.
Let $B$ denote the set of $\sqrt{\varepsilon}$-bad vertices of $H$.
First we want to match all the bad vertices.
Since $H$ is $\varepsilon$-close to $H_{1,3}'(n)$, we have $|B|\leq 4\sqrt{\varepsilon}n$. Let $U_i^{bad}=U_i\cap B$, $W_i^{bad}=W_i\cap B$ for $i\in [4]$.
Write $b':=\max\{|U_1^{bad}|, \sum_{i=2}^4|W_i^{bad}|\}$
and let $b\in \{b',b'+1,b'+2\}$ such that $n-b\equiv 0\pmod 3$.
Then we see that $b\equiv s\pmod 3$. Write $b=3t+s$.
Let $A\subseteq V(H)$ such that $(U_2\cup U_3\cup U_4\cup B)\subseteq A$, $|A\cap U_1|=b$, $|W_2\setminus A|=|W_3\setminus A|=|W_4\setminus A|=(n-b)/3$.
Let $H':=H[A]$.
Then for any $x\in A\cap V_1$,
\[
\delta_1(N_{H'}(x))>\delta((2n+b)/3,t,s)=\left\{\begin{array}{ll} (2n+b)^2/9-((2n+b)/3-t)^2 &\text{if}\
s=1 , \\[5pt]
 (2n+b)^2/9-((2n+b)/3-t+1)((2n+b)/3-t-1) &\text{if}\
s=2,\\[5pt]
 (2n+b)^2/9 - ((2n+b)/3-t)((2n+b)/3-t-1) &\text{if}\ s=3.
\end{array}\right.
\]
Then by Theorem \ref{small-matchings}, $\{N_{H'}(x)\ |\ x\in A\cap V_1\}$ has a rainbow matching. So $H'$ has a matching of size $b$ denoted by $M_1$.

Write $W_i^1=W_i\backslash(B\cup V(M))$ for $i=2,3,4$.
We choose $\eta$ such that $\varepsilon\ll \eta\ll 1$.
Write $B_1=B \setminus V(M_1)$. 
A vertex $v\in B_1$ is \textit{useful} if the number of edges containing $v$ and exactly one vertex in $(W_2^1\cup W_3^1\cup W_4^1)$ is at least $\eta n^3$.  We denote the set of useful vertices by $B_{11}$ and let $B_{12}=B_1 \setminus B_{11}$. Recall  that $|B_{11}|\leq 4\sqrt{\varepsilon}n$ and $\varepsilon\ll \eta$.
So we can greedily find a matching $M_{21}$ such that for any  $e\in M_{21}$ is type $UUUW$ and $|e\cap B_{11}|=1$.
Note that for $x\in B_{12}$, $d_{H-V(M_{1}\cup M_{21})}(x)>n(2nr-r^2)-3|M_{1}\cup M_{21}|n^2>4n^3/9$ and the number of edges containing $x$ and at least two vertices of
 $(W_2^1\cup W_3^1\cup W_4^1)$ is at most $3n(n/3)(n/3) = n^3/3$. So there are at least $n^3/9$ edges in $H-V(M_{1}\cup M_{21})-(W_2^1\cup W_3^1\cup W_4^1)$ containing $x$.
 Thus we may greedily find a matching $M_{22}$ of size $|B_{12}|$ in $H-V(M_1\cup M_{21})-(W_2^1\cup W_3^1\cup W_4^1)$ covering $B_{12}$.

Write $H_1=H-V(M_1 \cup M_{2})$, where $M_2=M_{21}\cup M_{22}$.
 Now every vertex of $V(H_1)$ in $H$ is $\sqrt{\varepsilon}$-good with respect to $H_{1,3}'(n)$.
 Next we match some vertices such that the ratio of the number of unmatched vertices in $W_2 \cup W_3 \cup W_4$ and those in $U_2 \cup U_3 \cup U_4$ is $1/3$.
 So we can greedily find a matching $M_3$ of size $|B_{12}|$ such that every edge in $M_{3}$ is type $UUWW$. For $2\leq i\leq 4$, $W_{i}^2=W_{i}^1\backslash V(M_{2}\cup M_3)$. One can see that $\sum_{i=2}^4|W_i^2|=n-|M_1\cup M_2\cup M_{3}|$.

Now we show that we can find a matching such that the number of unmatched vertices in $W_i$ is the same for $i=2,3,4$.
Let $H_2=H-V(M_3)$.  Without loss generality, suppose that $|W_2^2|\leq |W_3^2|\leq |W_4^2|$. It is easy to see  that $|W_3^2|-|W_2^2|\leq |W_4^2|-|W_2^2|\leq 3b$. Since every vertex of $W_2^2\cup W_{3}^2 \cup W_4^2$ in $H$ is $\sqrt{\varepsilon}$-good with respect to $H_{1,3}'(n)$, then for every $x\in W_{3}^2$, the number of edges containing $x$ in $H_2$ and avoiding $W_2 \cup W_4$ is at least
\[
4n^3/9-{\varepsilon}^{1/4}n^3-3|M_3\cup M_1\cup M_{21}\cup M_{22}|>4n^3/9-{\varepsilon}^{1/4}n^3-4bn^3>n^3/3.
\]
Thus we can find a matching $M_{31}$ in $H_2$ of size $|W_3^2|-|W_2^2|$ such every edges in $M_{31}$ containing exactly one vertex from $W_3^2$ and avoiding $W_2^2 \cup W_4^2$.
Similarly. we can also find a matching $M_{32}$ in  $H_2-V(M_{31})$ such that every edges in $M_{32}$ of size $|W_4^2|-|W_2^2|$ containing exactly one vertex from $W_4^2$ and avoiding $W_2^2 \cup W_3^2$.
Let $M_3:=M_{31}\cup M_{32}$ and for $2\leq i\leq 4$, $W_i^3:=W_i^2 \backslash V(M_3)$. Let $H_3:=H_2-V(M_3)$.

One can see that $|W_2^3|=|W_3^3|=|W_4^3|$ and $3|W_2^3|=|V(H_3)\cap V_1|$. Moreover, every vertex in $H_3$ is $\varepsilon^{1/5}$-good with $H_{1,3}'(n')$, where $n'=n-|\cup_{i=1}^3 M_i|$.
By Lemma \ref{good-PM}, $H_3$ has a perfect matching, say $M_4$.
Now $\cup_{i=1}^4M_i$ is a perfect matching of $H$. This completes the proof. \qed

\section{Absorbing lemma}

In this section, we prove an absorbing lemma for balanced $4$-graphs.
We need the following absorbing lemma from \cite{LM}.
\begin{lemma}[Lo and Markstr\"om, \cite{LM}]\label{absorbing lemma}
Let $0 <\tau < 1/(10k^3)$ and $\tau'= \tau^{2k-1}/20$. Then there is
an integer $n_0$ such that for all $n> n_0$  the following holds:
suppose $H$ is a $k$-partite $k$-graph with $n$ vertices in each
class  and minimum degree $\delta_1(H)\geq(1/2 +\tau)n^{k-1}$, then
there exists a matching $M$ in $H$ of size $|M|\leq(k-1)\tau^kn$
such that, for every balanced set $W$ of size $|W|\leq k\tau'n$,
there exists a matching covering exactly the vertices of $V(M)\cup
W$.
\end{lemma}

We adapt Lemma \ref{absorbing lemma} to balanced $4$-graphs and show that there exists a matching in balanced $4$-graphs which can absorb any balanced unmatched set of vertices of small size given degree conditions.

\begin{lemma}\label{absorbing lemma2}
Let $0<1/n\ll\beta\ll \gamma\ll 1$.
Let $H=(V_1,V_2,V_3,V_4)$ be a balanced 4-graph with $n$ vertices in each class such that $d_{H}(\{x,y\})\geq (1/2+\gamma)n^2$ for any $x\in V_1$ and $y\in V_2\cup V_3\cup V_4$.
Then $H$ have  a matching $M$ with $|M|\leq \gamma^6n$ such that for any balanced subset $S$ with $|S|\leq \beta n$, $H[S\cup V(M)]$ has a perfect matching.

\end{lemma}

\pf Since $d_{H}(\{x,y\})\geq (1/2+\gamma)n^2$ for any $x\in V_1$ and $y\in V_2\cup V_3\cup V_4$, then $d_H(v)\geq (1/2+\gamma)n^3$ for all vertex $v \in V(H)$, i.e., $\delta_1(H)\geq (1/2+\gamma)n^2$. Thus by Lemma \ref{absorbing lemma}, the statement follows. \qed

\section{Almost perfect matchings in non-extremal case}

In this section, we show that $H$ contains an almost perfect matching in the non-extremal case.
We sketch the proof as follows.
By a theorem of Pippenger (Theorem \ref{nibble-new}) on edge cover, it suffices to find an almost regular spanning subgraph $F$ of $H$ with small codegree.
Theorem \ref{AHJ19} allows us to find a fractional perfect matching with small support in $H$ (Lemma \ref{FPM}).
In order to obtain $F$, we construct many edge-disjoint fractional perfect matchings with small codegree (Lemma 5.5), and sample each edge with certain probability related to the values of those matchings. Now such $F$ has desired properties (Lemma 5.6).

First we prove that if $H$ is not close to the extremal graph, then every large set of vertices induce many edges.

\begin{lemma}\label{dense}
Let  $0<1/n\ll \gamma \ll  \varepsilon\ll 1$.
Let $H$ be a balanced $4$-partite $4$-graph on vertex set $V_1, V_2, V_3, V_4$ with $n$ vertices in each class such that $H$ is not $\varepsilon$-close to $H_{1,3}'(n)$ and   $d_{H}(\{x,y\})\geq (5/9-\gamma)n^2$ for any $x\in V_1$ and $y\in V_2\cup V_3\cup V_4$. Then for any vertex subset $U \subseteq V(H)$ such that $V_1\subseteq U$ and $|U\cap V_i|\geq 2n/3$ for $2\leq i\leq 4$, $e(H[U])\geq \varepsilon n^4/2.$
\end{lemma}

\pf Suppose that the result does dot hold. Then there exists $U$ with $V_1\subseteq U$ and $|U\cap V_i|= 2n/3$ for $2\leq i\leq 4$, such
that $e(H[U])\leq \varepsilon n^4 / 2$. The number of edges intersecting both $V(H) \setminus U$ and $U \setminus V_1$ is at least
\begin{align*}
&\frac{1}{2}\sum_{x\in V_1,y\in U-V_1}d_{H}(x,y,V(H)-U,U)+\sum_{x\in V_1,y\in U-V_1}d_{H}(x,y,V(H)-U,V(H)-U)-e(H[U])\\
&\geq \frac{1}{2}n(|U|-|V_1|)\frac{4}{9}n^2+n|U\cap V_2|(\frac{5}{9}n^2-\gamma n^2-|V_3-U||U\cap V_4|-|V_4-U||U\cap V_3|)\\
&+n|U\cap V_3|(\frac{5}{9}n^2-\gamma n^2-|V_2-U||U\cap V_4|-|V_4-U||U\cap V_2|)\\
&+n|U\cap V_4|(\frac{5}{9}n^2-\gamma n^2-|V_2-U||U\cap V_3|-|V_3-U||U\cap V_2|)-\varepsilon n^4/2\\
&=\frac{1}{2}n(|U|-|V_1|)\frac{4}{9}n^2+n|U-V_1|(\frac{5}{9}n^2-2|V_2-U||U\cap V_4|)\\
&=\frac{4}{9}n^4+2n^2(\frac{1}{9}n^2-\gamma n^2)-\varepsilon n^4/2\\
&=\frac{4}{9}n^4+\frac{2}{9}n^3-2\gamma n^4-\varepsilon n^4/2\\
&> \frac{2}{3}n^4-\varepsilon n^4,
\end{align*}
which implies that $H$ is $\varepsilon$-close to $H_{1,3}'(n)$, a contradiction. \qed

\begin{comment}
We also need the following definition.

\begin{definition}
Let $n,k,l \ge 2$ be three integers.
Let $d_1, \ldots, d_k$ be integers such that $0 \le d_i \le n$ for $i \in [k]$.
Let $V_1, \ldots, V_k$ be $k$ sets of size $n$ and $V = V_1 \cup \ldots \cup V_k$.
Let $W_i \subseteq V_i$ with $|W_i| = d_i$ for each $i \in [k]$ and  $W = W_1 \cup \ldots \cup W_k$.
Let $H'_{k,l}(n;d_1, \ldots, d_k)$ be $k$-partite $l$-graph with vertex partition $V_1, \ldots, V_k$ and
edge set
$$E(H'_{k,l}(n;d_1, \ldots, d_k) = \{e \in V_{i_1} \times \ldots \times V_{i_k}   | \{i_1,\ldots,i_k\} \in {[k] \choose l}, e \cap W \ne \emptyset, e \cap (V \setminus W) \ne \emptyset \}.$$
We also denote $H'_{k,k}(n;d_1, \ldots, d_k)$ by $H'_{k}(n;d_1, \ldots, d_k)$.
\end{definition}

\end{comment}

%%%%%%%%%%%%%%%%%%%%%%%%%%%%%%%%%%%%%%%%%%%%%%%%%%%%%%%%%%%%%%%%%%%

Let $V_i:=\{v_{ij}\ |\ j\in [n]\}$ for $i\in [r]$. Let $H$ be an $r$-graph with vertex set $\cup_{i=1}^r V_{i}$. $H$ is \textit{stable} if  $j_s\leq k_s$ for $s\in [r]$, $\{u_{1k_1},\ldots, u_{rk_r}\}\in E(H)$ implies that $\{u_{1j_1},\ldots, u_{rj_r}\}\in E(H)$.
Now we establish a stability result of $3$-partite stable graphs.

\begin{lemma}\label{frac-stab}
Let $G$ be a  3-partite 2-graph on vertex set $V_1,V_2,V_3$ such that $|V_i|=n$ for $i \in [3]$ and $G[V_i\cup V_j]$ is stable for $\{i,j\}\in {[3]\choose 2}$. If $e(G[V_i\cup V_j])\geq (5/9-\gamma)n^2$ for any $\{i,j\}\in {[3]\choose 2}$ and $\nu(G)\leq n-1$, then $G$ is $16\sqrt{\gamma}$-close to $H_{3,2}(n;n/3,n/3,n/3)$.
\end{lemma}

% Tutte-Berge formula: max matching size = 1/2 * min_{U \subseteq V} (|U| - odd(G-U) + |V|)

\pf Since $\nu(G)\leq n-1$, by Tutte-Berge Formula, there exists $S\subseteq V(G)$ such that $c(G-S)-|S|\geq n+ 2$ where $c(G-S)$ is the number of odd component in $G-S$.
We choose $S$ to be maximal among all the choices of $S$.
%We choose $S$ to be maximal.
%Since $G$ is stable, we have the following claim.

\medskip
\textbf{Claim~1.} $G-S$ contains at most one non-trivial component.

Otherwise, suppose that $G-S$ has two non-trivial components say $C_1, C_2$. We claim that $C_1$ and $C_2$ are non-bipartite; For otherwise, we could move all the vertices in the smaller part of $C_1$ (or $C_2$) to $S$, which would contradict the maximality of $S$.
Let $u_{i_1j_1}u_{i_2j_2}\in E(C_1)$ and $u_{i_1k_1}u_{i_2k_2}\in E(C_2)$. If $j_1\leq k_1$, then since $G[V_i \cup V_j]$ is stable for $\{i,j\} \in {3 \choose 2}$, $u_{i_1j_1}u_{i_2k_2}\in E(G)$, a contradiction as $C_1, C_2$ are disconnected.
The other case can be treated similarly.
%Else, $u_{i_1k_1}u_{i_2j_2}\in E(G)$, a contradiction again.

\medskip

Let $C$ be the non-trivial component in $G-S$.
Write $s:=|S|$, $c:=|V(C)|$, $s_i:=|S\cap V_i|$ and $c_i:=|V(C)\cap V_i|$. By Claim 1, one can see that 
\[
\nu(G) \le s+(c-1)/2\leq n-1,
 \]which implies that $(c-1)/2\leq n-1-s$.  Without loss generality, suppose that $s_1+c_1/2\geq s_2+c_2/2\geq s_3+c_3/2$. Recall that $C$ is factor-critical. Thus we have $c/2>c_i>0$ for $i\in [3]$. %So we may conclude
%%\[
%r_2+m_3\leq (r-1)/2+m_3\leq n-1-m_1-m_2,
%\]
%and
%\[
%r_3+m_2\leq (r-1)/2+m_3\leq n-1-m_1-m_3.
%\]
 Let $G_3:=G[V_2\cup V_3]$. Note that $s_2+c_2/2\leq (n-s_3-c_3/2)/2$.  Then
 \begin{align*}
(\frac{5}{9}-\gamma)n^2< e(G_3)&\leq c_2c_3+n^2-(n-s_2)(n-s_3) \text{\quad (by Claim 1)}  \\
&=c_2c_3+(s_2+s_3)n-s_2s_3 \\
&\leq c_2c_3+(s_3+(n-s_3-c_3/2)/2-c_2/2)n-((n-s_3-c_3/2)/2-c_2/2)s_3\\
&=n^2/2+s_3^2/2-(c_3/4+c_2/2)(n-s_3)+c_2c_3\\
&= n^2/2+s_3^2/2-(c_3/2+c_2)(n-s_3)/2+c_2c_3\\
&\leq n^2/2+s_3^2/2, \text{\quad (as $n-s_3 \ge c/2 > c_3$)}
\end{align*}
 which implies that $s_3\geq (1/3-\sqrt{\gamma})n$. Recall that $s_3+c_3/2\leq n/3 - 1/6$. So we may infer that $c_3\leq 2\sqrt{\gamma} n$. With similar discussion, we have $s_2\geq (1/3-\sqrt{\gamma})n$ and $c_2\leq 3\sqrt{\gamma} n$. Since  $C$ is factor-critical, $c_1\leq c_2+c_3\leq 5\sqrt{\gamma} n$. Note that $s_1+c_1/2\geq (1/3-\sqrt{\gamma})n$. Hence we have $s_1\geq (1/3-7\sqrt{\gamma}/2)n$.

Now we have
\[
|H_{3,2}(n;n/3,n/3,n/3)\backslash E(G)|\leq \frac{11}{2}\sqrt{\gamma} n^2+(c_1+c_2+c_3)n^2\leq 16 \sqrt{\gamma} n^2.
\]
Thus $G$ is $16 \sqrt{\gamma}$-close to $H_{3,2}(n; n/3,n/3,n/3)$. This completes the proof. \qed

%We aim to construct such an almost perfect matching from many fractional perfect matchings.

A {\it fractional matching} in a $k$-graph
  $H$ is a function $f: E\rightarrow [0,1]$ such that for any
  $v\in V(H)$, $\sum_{\{e\in E: v\in e\}}f(e)\le 1$. The size of $f$ is
$\sum_{e\in E(H)}f(e)$. A {\it fractional matching} $f$ is called \emph{maximum fractional matching} if $\sum_{e\in E(H)}f(e)\geq \sum_{e\in E(H)}f'(e)$ for any fractional matching $f'$. Let $\nu_f(H)$ denote the size of a maximum fractional matching.   A fractional matching $f$ is  {\it perfect} if $\sum_{e\in E}f(e)=|V(H)|/k$.
   A {\it fractional vertex cover} of $H$ is a function   $\omega:V(H)\rightarrow \mathbb{R}^+\cup \{0\}$ such that $\sum_{x\in e}\omega(x)\ge 1$
for $e\in E(H)$. A fractional vertex cover  $\omega$ of $H$ is called \emph{minimum fractional vertex cover} if for any fractional vertex cover $\omega'$, $\sum_{v\in V(H)}\omega(v)\leq \sum_{v\in V(H)}\omega'(v)$. We use $\mu_f(H)$ to denote the size of a minimum fractional vertex cover.  By linear programming duality theory, one can see that $\nu_f(H)=\mu_f(H)$ for any hypergraph $H$.

  Given a family $\{F_1,\ldots, F_m\}$  of hypergraphs, a \emph{rainbow fractional matching} of size $n$ of $\{F_1,\ldots, F_m\}$ is a set of edges $e_1\in E(F_1),\ldots, e_m\in E(F_m)$ together with a fractional matching $f:\{e_1,\ldots,e_m\}\rightarrow [0,1]$ of size $n$.  A  rainbow fractional matching $f$ is called a \emph{rainbow perfect fractional matching} if $V(F_1)=\cdots=V(F_m)$ and $\sum_{v\in e}f(e)=1$ for any $v\in V(F_1)$.

We need the following theorem on rainbow fractional matchings in \cite{AHJ}.
\begin{theorem}[Aharoni, Holzman, Jiang, Theorem 1.6 in \cite{AHJ}]\label{AHJ19}
Let $r \geq 2$ be an integer, and let $n$ be a positive rational number. Let $H_1,\ldots, H_{\lceil rn \rceil }$ be $r$-graphs such that $\nu_f(H_i)\geq n$ for $i= 1,...,\lceil rn \rceil $.
Then $H_1,\ldots, H_{\lceil rn \rceil }$ has a rainbow fractional perfect matching of size $n$.
\end{theorem}

Using Theorem \ref{AHJ19}, we manage to show that $H$ contains a fractional perfect matching whose support is small.

\begin{lemma}\label{FPM}
Let  $0<1/n\ll \gamma\ll \varepsilon\ll 1$.
Let $H$ be a balanced $4$-partite $4$-graph on vertex set $V_1,V_2,V_3,V_4$ with $n$ vertices in each class such that $H$ is not $\varepsilon$-close to $H_{1,3}'(n)$ and   $d_{H}(\{x,y\})\geq (5/9-\gamma)n^2$ for any $x\in V_1$ and $y\in V_2\cup V_3\cup V_4$. Then $H$ has a fractional perfect matching $f$ such that $|\{e\ |\ f(e)>0\} |\leq 4n$.
\end{lemma}
\pf By Theorem \ref{AHJ19} with $r=4$ and $H_1= \ldots = H_{4n}=H$, it suffices to show that $\nu_f(H)=n$.
Let
$\omega:V(H)\rightarrow [0,1]$ be a minimum fractional cover of $H$ such that for $i\in [4]$, $\omega(v_{i1})\geq \cdots\geq \omega(v_{in})$. Let
$Cl(H)$ be a balanced $4$-partite $4$-graph with vertex set $V(H)$ and edge set
\[
E(Cl(H)):=\{S\in V_1\times V_2\times V_3\times V_4\ |\ \sum_{v\in S}\omega(v)\geq 1\}.
\]
It is easy to see that $H$ is a subgraph of $Cl(H)$.
We aim to show that $\nu(Cl(H))=n$.

Let $G$ be a $3$-partite $2$-graph with vertex set $V_2\cup V_3\cup V_4$ and edge set $N_{Cl(H)}(\{v_{1n},v_{2n}\})\cup N_{Cl(H)}(\{v_{1n},v_{3n}\})\cup N_{Cl(H)}(\{v_{1n},v_{4n}\})$. Next we discuss two cases.

\medskip
\textbf{Case 1.~} $\nu(G)\geq n$.
\medskip

Let $M$ be a matching of $G$ of size $n$. For $2\leq i\leq 4$, let $M_i:=M\cap N_{Cl(H)}(\{v_{1n},v_{in}\})$ and  $m_i:=|M_i|$. Now we partition $V(H) \setminus V(M)$ into $n$ balanced 2-sets, say $A^2_1,\ldots, A^2_{m_2}, A^3_{1},\ldots,$ $A^3_{m_3},$ $A^4_1\ldots, A^4_{m_4}$ such that $A_1^i,\ldots,A_{m_i}^i\subseteq  (V_1\times V_i)$ for $2\leq i\leq 4$. By the definition of $Cl(H)$, one can see that $M_i\subseteq N_{H}(A_{j}^i)$ for $1\leq j\leq m_i$ and $2\leq i\leq 4$. For $2\leq i\leq 4$, write $M_i=\{e_1^i,\ldots,e_{m_i}^i\}$. Then
\[
\bigcup_{i=2}^4\{A_j^i\cup e_j^i\ |\ 1\leq j\leq m_i\}
\]
is a perfect matching of $H$.

\medskip
\textbf{Case 2.~$\nu(G)< n$}
\medskip

By Lemma \ref{frac-stab}, $G$ is $16 \sqrt{\gamma}$-close to $H_{3,2}(n;n/3,n/3,n/3)$. Then there exists $T\subseteq V(H)-V_1$ with $|T\cap V_i| = (1/3-{\gamma}^{1/4})n$ for  $i\in \{2,3,4\}$ such that  for every $v\in T$, $v$ is $16\gamma^{1/4}$-good with respect to $H_{3,2}(n;n/3,n/3,n/3)$, which implies that $d_{G}(v)\geq (2 - 2\gamma^{1/4})n$ for all $v \in T$.

Since $H$ is not $\varepsilon$-close to $H_{1,3}'(n)$ and $d_{H}(\{x,y\})\geq (5/9-\gamma)n^2$ for any $x\in V_1$ and $y\in V_2\cup V_3\cup V_4$, by Lemma \ref{dense}, $H-T$ contains at least $\varepsilon n^4/2$ edges. So we can greedily find a matching $M_1$ of size $3{\gamma}^{1/4}n$ from $H-T$ as $\gamma \ll \varepsilon$.

Write $T :=\{x_1,\ldots,x_{n-3{\gamma}^{1/4}n}\}$. Next we can greedily find a matching $M_2'$ of size $n-3{\gamma}^{1/4}n$ in $G - V(M_1)$ such that every edge of $M_2'$ intersecting $T$ exactly one vertex since $d_{G}(v)\geq (2 - 2\gamma^{1/4})n$ for all $v \in T$.
Write $M_2'=\{e_1,\ldots,e_{n-3{\gamma}^{1/4}n}\}$.

Note that we may partition $V(H)-V(M_1)-V(M_2')$ into $n-3{\gamma}^{1/4}n$ 2-sets $A_1,\ldots, A_{n-3{\gamma}^{1/4}n}$ such that $|A_i\cap V_1|=1$. By the definition of $Cl(H)$,
\[
M_2=\{A_i\cup e_i\ |\ i\in \{1,\ldots,n-3{\gamma}^{1/4}n\}
\]
is a matching of size $n-3{\gamma}^{1/4}n$ of $Cl(H)$.
One can see that $M_1\cup M_2$ is a perfect matching of $Cl(H)$.

In both cases, $Cl(H)$ has a perfect matching.
Thus $H$ has a fractional perfect matching. Let $H_i:=H$ for $1\leq i\leq n$. Applying Theorem \ref{AHJ19} to $H_i$, $H$ has a fractional perfect matching $f$ such that $|\{e\in E(H)\ |\ f(e)>0\}|\leq 4n$.
This completes the proof. \qed

Next we construct many edge-disjoint fractional perfect matchings such that the codegree of any two vertices is small.

\begin{lemma}\label{dis-FPM}
Let  $0<1/n\ll \gamma\ll \varepsilon\ll 1$.
Let $H$ be a balanced $4$-partite $4$-graph on vertex set $V_1,V_2,V_3,V_4$ with $n$ vertices in each class such that $H$ is not $\varepsilon$-close to $H_{1,3}'(n)$ and   $d_{H}(\{x,y\})\geq (5/9-\gamma)n^2$ for any $x\in V_1$ and $y\in V_2\cup V_3\cup V_4$. Then $H$ has $n/\ln n$ edge-disjoint fractional perfect matchings $f_1,\ldots,f_{n/\ln n}$ such that
for any two vertices   $x,y \in V(H)$, $\sum_{i=1}^{n/\ln n}\sum_{\{x,y\}\subseteq e}f_i(e)\leq 3$.
%and for any $v\in V(H)$, $\sum_{i=1}^{n/\ln n}\sum_{v\in e}f_{i}(e)=n/\ln n$.
\end{lemma}

\pf By Lemma \ref{FPM}, there exists a fractional perfect matching $f_1$ in $H$ such that $|\{e\in E(H)\ |\ f_1(e)>0\}|\leq 4n.$
Write $M_1:=\{e\in E(H)\ |\ f_1(e)>0\}$.
%and $M_0':=\{e\in E(H)\ |\ f_0(e)\}$.
Let $s \ge 1$ be the maximum integer such that there exist $s$ fractional perfect matchings $f_1,\ldots, f_s$ such that for $i\in [s]$, $M_i:=\{e\in E(H)\ |\ f_i(e)>0\}$ with $|M_i|\leq 4n$ and   for any $\{i,j\}\in {[s]\choose 2}$, $M_i\cap M_j=\emptyset$ and for any $\{x,y\}\in {V(H)\choose 2}$, $\sum_{i=1}^s \sum_{\{x,y\}\subseteq e\in E(H)}f_i(e)\leq 3$.

By contradiction, suppose $r < n/\ln n$.
Let
\[
A_{s+1}:=\{\{x,y\}\in {V(H)\choose 2}\ |\ \sum_{i=1}^s \sum_{\{x,y\}\subseteq e\in E(H)}f_i(e)>2\},
\]
\[
B_{s+1}:=\{ e\in E(H_r)\ |\ \exists S\in A_{s+1}\ \mbox{such that } S\subseteq e\}
\]
%Put $H_{r+1}:=H_r-M_r-B_{r+1}$.
and
$$H_{s+1} := H \setminus E(M_1 \cup \cdots M_s \cup B_{s+1}).$$

\medskip
\textbf{Claim 1.}~$H_{s+1}$ has a fractional perfect matching.

Let $\omega_{s+1}$ be a fractional vertex cover such that such that for $i\in [4]$, $\omega_{s+1}(v_{i1})\geq\cdots\geq \omega_{s+1}(v_{in})$.
For $x\in V_1$, since $\sum_{i=1}^s\sum_{x\in e}f_i(e)=s$, there are at  most $s/2$ vertices $y\in V_2\cup V_3\cup V_4$ such that
$\sum_{i=1}^s \sum_{\{x,y\}\subseteq e\in E(H)}f_i(e)>2$. Let
\[
T:=\{y\in V(H)-V_1\ |\ \sum_{i=1}^s \sum_{\{v_{1n},y\}\subseteq e\in E(H)}f_i(e)>2\}.
\]
Write $t:=|T|$ and $T=\{u_1,\ldots,u_{t}\}$. So $t\leq s/2$. Let $Cl(H_{s+1})$ be a $4$-partite $4$-graphs with vertex set $V(H)$ and edge set
\[
E(Cl(H_{s+1}))=\{e\in \prod_{i=1}^4V_i\ |\ \sum_{x\in e}\omega_{s+1}(x)\geq 1\}.
\]

We aim to show that $Cl(H_{r+1})$ has a fractional perfect matching.
For any $x\in V(H) \setminus (V_1 \cup T)$, one can see that
\begin{align*}
d_{Cl(H_{s+1})}(\{v_{1n},x\})&\geq d_{H_{s+1}}(\{v_{1n},x\})\\
&\geq d_{H}(\{v_{1n},x\})-\sum_{i=1}^s|M_i|-|\{e\in E(H)\ |\ \{v_{1n},x\}\subseteq e \mbox{ and }e\cap T \ne \emptyset\}|\\
&\geq (5/9-\gamma)n^2-4n^2/\ln n-n^2/(2\ln n)\\
&\geq (5/9-2\gamma)n^2.
\end{align*}

With similar discussion, for any $y\in T$, there are at most $s/2$ vertices $x$ such that $\sum_{i=1}^s \sum_{\{x,y\}\subseteq e\in E(H)}f_i(e)>2$.
Thus
\begin{align*}
d_{Cl(H_{s+1})}(y)&\geq d_{H_{s+1}}(y)\\
&\geq d_{H}(y)-\sum_{i=1}^s|M_i|-n^2 (s/2)\\
&\geq\frac{1}{3} (5/9-\gamma)n^3-4n^2/\ln n-n^3/(2\ln n)\\
&\geq \frac{1}{3}(5/9-2\gamma)n^3.
\end{align*}
Hence we may greedily find a matching $\mathcal{M}$ of size $t$ covering $T$. By stability, we may assume that $v_{1n}\notin V(\mathcal{M}) - V_1$.  One can see that for any $x\in V(H)-V(\mathcal{M})-V_1$,
\begin{align*}
d_{Cl(H_{s+1})-V(\mathcal{M})}(\{v_{1n},x\})&\geq d_{Cl(H_{s+1})-V(\mathcal{M})}-3|\mathcal{M}|n\\
&\geq (5/9-3\gamma)n^2.
\end{align*}
Moreover,
since $H_{s+1}$ is not $\varepsilon$-close to $H_4(n;n/3,n/3,n/3)$,
$\mathcal{H}_{s+1}-V(\mathcal{M})$ is not $\frac{\varepsilon}{2}$-close to  $H_4(n-|T|,n/3-|T|/3,n/3-|T|/3,n/3-|T|/3).$ By Lemma \ref{FPM}, $\mathcal{H}_{s+1}-V(\mathcal{M})$ has a fractional perfect matching $f_{s+1}'$.
So $\mathcal{M} \cup \{f_{s+1}'\}$ is a fractional perfect matching of $Cl(\mathcal{H}_{s+1})$.
By linear programming duality,
 $\mathcal{H}_{s+1}$ has a fractional perfect matching.

\medskip

By Theorem \ref{AHJ19}, there exists a fractional perfect matching $f_{s+1}$ in $H_{s+1}$ such that $|\{e\in E(H_{s+1})\ |\ f_{s+1}(e)>0\}|\leq 4n$.
By construction, $f_1,\ldots, f_{s+1}$  contradicts to the maximality of $s$. \qed
%Continuing the process, we may find $n/\ln n$ fractional perfect matchings satisfying the given conditions. \qed

By sampling the edge-disjoint fractional perfect matchings of $H$ obtained from Lemma \ref{dis-FPM}, we have an almost regular spanning subgraph $F$ of $H$ with small codegree.

\begin{lemma} \label{lem-sampling}
Let  $0<1/n\ll \gamma\ll \varepsilon\ll 1$.
Let $H=(V_1,V_2,V_3,V_4)$ be a balanced $4$-partite $4$-graph with $n$ vertices in each class such that $H$ is not $\varepsilon$-close to $H_{1,3}'(n)$ and   $d_{H}(\{x,y\})\geq (5/9-\gamma)n^2$ for any $x\in V_1$ and $y\in V_2\cup V_3\cup V_4$. Then there exists a spanning subgraph $F$ such that
         \begin{itemize}
            \item [$(i)$] $d_{F}(v) =(1+o(1))n/\ln n$ for $v\in V(H)$, and
            \item [$(ii)$] $\Delta_2(F)\le n^{0.1}$.
         \end{itemize}
\end{lemma}

\pf By Lemma \ref{dis-FPM}, $H$ has $n/\ln n$ edge-disjoint fractional perfect matching $f_1,\ldots,f_{n/\ln n}$ such that for any $\{x,y\}\in {V(H)\choose 2}$,
\[
\sum_{i=1}^{n/\ln n}\sum_{\{x,y\}\subseteq e}f_i(e) \le 3.
\]
Note that $\sum_{i=1}^{n/\ln n}f_i(e)\leq 1$ for any $e\in E(H)$.
Let $F$ be a spanning subgraph of $H$ obtained by
        independently selecting each edge $e$ at random
    with probability $\sum_{i=1}^{n/\ln n}f_i(e)$. Hence, for $v\in V(H)$ and  $\{x,y\}\in {V(H)\choose 2}$,
	\[
		d_{F}(v)=\sum _{v\in e} X_e \mbox{ and } d_{F}(\{x,y\})= \sum _{\{x,y\}\subseteq e} X_e
	\]
	where $X_e\sim Ber(\sum_{i=1}^{n/\ln n}f_i(e))$ is the Bernoulli random variable with $X_e=1$ if $e\in E(F)$ and $X_e=0$ otherwise.
 Thus,
    since $f_1,\ldots, f_{n/\ln n}$ are $n/\ln n$  are edge-disjoint perfect fractional  matching in $H$, for any $v\in V(H)$
	\[
		\E (d_{F}(v))=\sum_{v\in e} \mathbb{P}(X_e) = \sum_{i=1}^{n/\ln n} \sum _{v\in e} f_i(e)=n/\ln n,
	\]
and for any $\{x,y\}\in {V(H)\choose 2}$,
\[
		\E (d_{F}(\{x,y\}))=\sum_{\{x,y\}\subseteq e} \mathbb{P}(X_e)=\sum_{i=1}^{n/\ln n} \sum _{\{x,y\}\subseteq e} f_i(e)\leq 3.
	\]

Now by the Chernoff bound, for any $v\in V(H)$ and  $\{x,y\}\in {V(H)\choose 2}$,
\[
		\P(|d_{F}(v)-n/\ln n|\ge n^{0.9})\le e^{-\Omega(n^{0.5})},
	\]
and
\[
		\P(d_{F}(\{x,y\})\ge n^{0.1})\le e^{-\Omega(n^{0.1})}.
	\]
Hence by union bound, we have $\Delta_2({F})\le n^{0.1}$ and  $n/\ln n-n^{0.9}\leq d_{F}(v)\leq n/\ln n+n^{0.9}$ with probability $1-o(1)$.
Therefore, with probability $1-o(1)$, ${F}$ satisfies $(i)$ and $(ii)$. This completes the proof. \qed

We also need a theorem due to Pippenger on the edge cover of hypergraphs.

\begin{theorem}[Pippenger \cite{PS}]
\label{nibble-new}
	For every integer $k\ge 2$ and reals $r\ge 1$ and $a>0$, there are $\gamma=\gamma(k,r,a)>0$ and $d_0=d_0(k,r,a)$ such that for every $n$ and $D\ge d_0$ the following holds: Every $k$-uniform hypergraph $H=(V,E)$ on a set $V$ of $n$ vertices in which all vertices have positive degrees and which satisfies the following conditions:
	\begin{itemize}
	 	\setlength{\itemsep}{0pt}
		\setlength{\parsep}{0pt}
		\setlength{\parskip}{0pt}
		\item[$(1)$] for all vertices $x\in V$ but at most $\gamma n$ of them, $d_H(x)=(1\pm \gamma)D$;
		\item[$(2)$] for all $x\in V$, $d_H(x)<r D$;
		\item[$(3)$] for any two distinct $x,y\in V$, $d_H(\{x,y\})<\gamma D$;
	\end{itemize}
	contains an edge cover of at most $(1+a)(n/k)$ edges.
\end{theorem}

Finally we can prove the main theorem of this section which shows that $H$ has an almost perfect matching when it is not close to the extremal graph.

\begin{theorem}\label{non-close-PM}
Let  $n,r,s$ be nonnegative integers such that $0<1/n\ll\gamma \ll \varepsilon\ll 1$.
Let $H$ be
$4$-partite $4$-graphs with $n$ vertices in each partition class and with partition  $V_1\cup V_2\cup V_3\cup V_4$.
If $H$ is not $\varepsilon$-close to $H_{1,3}'(n)$  and $d_{H}(\{x,y\}) >(5/9-\gamma)n^2$ for  $x\in V_1$ and $y\in V_2\cup V_3\cup V_4$, then $H$ has a matching of size at least $(1 - \eta)n$.
\end{theorem}

\pf
By Lemma \ref{lem-sampling}, there exists a spanning subgraph $F$ such that (i) $d_{F}(v) =(1+o(1))n/\ln n$ for $v\in V(H)$, and (ii) $\Delta_2(F)\le n^{0.1}$.
Applying Lemma \ref{nibble-new} to $F$, we obtain an edge cover $C$ of at most $(1+\eta/4)n$ edges.
Let $V'$ be the set of vertices that are incident to at least two edges in $C$.
Remove all the edges from $C$ that contain at least one vertex in $V'$, we obtain a matching of size at least $(1 - \eta)n$.
\qed

\section{Proof of Theorem \ref{main}}

\noindent \textbf{Proof of Theorem \ref{main}.}
By Observation \ref{Obser1}, it suffices to show that $H := H_{4}({\cal F})$ has a perfect matching.

Let $V_1,V_2,V_3,V_4$ be a partition corresponding to the definition of $H$.
Let $0<1/n \ll \eta \ll\beta\ll \gamma\ll \varepsilon \ll 1$.
Note that $d_{H}(\{x,y\}) >\delta(n,r,s)$ for  $x\in V_1$ and $y\in V_2\cup V_3\cup V_4$.
First suppose $H$ is $\varepsilon$-close to $H_{1,3}'(n)$.
By Theorem \ref{close-PM}, $H$ has a perfect matching.

So $H$ is not $\varepsilon$-close to $H_{1,3}'(n)$.
Note that $d_{H}(\{x,y\}) \ge 5n^2/9 \geq (1/2+\gamma)n^2$ for any $x\in V_1$ and $y\in V_2\cup V_3\cup V_4$.
By Lemma \ref{absorbing lemma2}, $H$ has  a matching $M_1$ with $|M_1|\leq \gamma^6n$ such that for any balanced subset $S$ with $|S|\leq \beta n$, $H[S\cup V(M_1)]$ has a perfect matching.

Write $H' := H - V(M_1)$ and $n':=|V_1-V(M_1)|$.
It is easy to see that $H'$ is not $\varepsilon/2$-close to $H_{1,3}'(n')$ and $d_{H}(\{x,y\}) >(5/9-\gamma)n^2$ for  $x\in V_1$ and $y\in V_2\cup V_3\cup V_4$.
Then by Theorem \ref{non-close-PM}, $H'$ has a matching $M_2$ of size at least $(1 - \eta)n$.
Write $S := V(H)' \setminus V(M_2)$.
By definition of $M_1$, $H[S \cup V(M_1)]$ a perfect matching $M_3$.
Therefore, $M_2 \cup M_3$ is a perfect matching of $H$.
This completes the proof of Theorem \ref{main}.
\qed

%\section{Appendix}

%
%
%
%{\color{blue}
%\begin{conjecture}
%Let $H$ be a $3$-partite 3-graph with partition $(V_1,V_2,V_3)$, where $|V_i|\in \{n,n+1,n+2\}$. If $\nu(H)=2$ and $\delta(H)\geq 3$, then there exists two vertices $x,y$ such that every edge of $H$ intersect with $\{x,y\}$.
%\end{conjecture}
%
%
%}


\begin{thebibliography}{99}

\bibitem{AGS} R. Aharoni, A. Georgakopoulos and P. Spr\"ussel, Perfect matchings in $r$-partite $r$-graphs, \emph{European J. Combin.}, \textbf{30} (2009), 39--42.

\bibitem{AHJ} R. Aharoni, R. Holzman and Z. Jiang, Rainbow fractional matchings, \emph{Combinatorica}, \textbf{39} (2019), 1191--1202.

%\bibitem{AH}  R. Aharoni and D. Howard, Size conditions for the
%existence of rainbow matching, preprint.

\bibitem{BDE} B. Bollob\'as, D.E. Daykin and P. Erd\"os, Sets of independent edges of a hypergraphs, \emph{Quart. J. Math. Oxford Ser.}, \textbf{27} (1976), 25--32.


\bibitem{CHWW2023} Y. Cheng, J. Han, B. Wang, and G. Wang, Rainbow spanning structures in graph and hypergraph systems. \textit{Forum of Mathematics, Sigma}, \textbf{11} (2023).

\bibitem{Diestel} R. Diestel, Graph Theory, Third edition, \textbf{173} (2005), Springer, Berlin.



\bibitem{Edmonds}J. Edmonds, Paths, trees and flowers, \emph{Canada. J. Math.}, \textbf{17} (1965) 449--467.

\bibitem{Erdos} P. Erd\H{o}s, On extremal problems of graphs and generalized
graphs, \emph{Israel J. Math.}, \textbf{2} (1964), 183--190.



\bibitem{HPS} H. H\'an, Y. Person and M. Schacht, On perfect matchings in
uniform hypergraphs with large minimum vertex degree, \emph{SIAM J.
Discret. Math.}, \textbf{23} (2009), 732--748.

\bibitem{Kan} I. Khan, Perfect matching in 3-uniform hypergraphs with large vertex degree,
\emph{SIAM. J. Discrete Math.}, \textbf{27} (2013), 1021--1039.


\bibitem{KST} P. K\"ov\'ari,
V. T. S\'os, P. Tur\'an, On a problem of Zarankiewicz, \emph{Colloq.
Math.}, \textbf{3} (1954), 50--57.


\bibitem{KO} D. K\"uhn and D. Osthus, Matchings in hypergraphs of large
minimum degree, \emph{J. Graph Theory}, \textbf{51} (2006),
269--280.

\bibitem{KOT} D. K\"uhn, D. Osthus and A. Treglown, Matchings in
3-uniform hypergraphs, \emph{J. Combin. Theory, Ser. B},
\textbf{103} (2013), 291--305.



\bibitem{LM} A. Lo and K. Markstr\"om, Perfect matchings in 3-partite 3-uniform hypergraphs,
\emph{J. Combin. Theory Ser. A}, \textbf{127} (2014), 22--57.

\bibitem{LWY2019} H. Lu, Y. Wang and X. Yu, Minimum codegree condition for perfect matchings in $k$-partite $k$-graphs, \textit{J. Graph Theory}, \textbf{92}(3) (2019), 207--229.


\bibitem{LWY2018} H. Lu, Y. Wang and X. Yu, Almost perfect matchings in $k$-partite $k$-graphs, \textit{SIAM J.
Discret. Math.}, \textbf{32}(1) (2018), 522--533.


\bibitem{LWY2022} H. Lu, Y. Wang and X. Yu, Rainbow perfect matchings for 4-uniform hypergraphs, \textit{SIAM J.
Discret. Math.}, \textbf{36}(3) (2022), 1645--1662.

\bibitem{LWY2023jcta} H. Lu, Y. Wang and X. Yu, A better bound on the size
  of rainbow matchings, \emph{J. Combin. Theory Ser. A},  \textbf{195} (2023), 105700.

\bibitem{LWY2023jctb} H. Lu, Y. Wang and X. Yu, Co-degree threshold for rainbow perfect matchings in uniform hypergraphs, \emph{J. Combin. Theory Ser. B}, \textbf{163} (2023), 83--111.

\bibitem{LWY2024} H. Lu, Y. Wang and X. Yu, On stability of rainbow matchings, arXiv:2302.06146, \textit{Submitted}.

\bibitem{LZ17} H. Lu and L. Zhang, Note on matchings in 3-partite 3-uniform hypergraphs, \emph{Discrete Math.}, \textbf{340} (5) (2017), 1136--1142.


\bibitem{Pikhurko} O. Pikhurko, Perfect
matchings and $K_4^{(3)}$-tilings in hypergraphs of large codegree,
\emph{Graphs Combin.}, \textbf{24} (2008), 391--404.

\bibitem{PS}
N. Pippenger and J. Spencer, Asymptotic behaviour of the chromatic index for hypergraphs, \emph{J. Combin. Theory, Ser. A}, \textbf{51} (1989), 24--42.


\bibitem{RR} V. R\"odl and A.
Ruci\'nski, Dirac-type questions for hypergraphs: a survey (or more
problems for Endre to solve), An Irregular Mind (Szemer\'edi is 70),
\emph{Bolyai Soc. Math. Studies}, \textbf{21} (2010).


\bibitem{RRS} V. R\"odl, A.
Ruci\'nski and E. Szemer\'edi, A Dirac-type theorem for 3-uniform
hypergraphs, \emph{Comb. Probab. Comput.}, \textbf{15} (2006),
229--251.

\bibitem{RRE} V. R\"odl, A. Ruci\'nski and E. Szemer\'edi, Perfect matchings in
uniform hypergraphs with large minimum degree, \emph{Europ. J.
Combin.}, \textbf{27} (2006), 1333--1349.



\bibitem{RRS2} V. R\"odl, A. Ruci\'nski and E. Szemer\'edi, An approximate
Dirac-type theorem for $k$-uniform hypergraphs,
\emph{Combinatorica}, \textbf{28 }(2008), 229--260.

\bibitem{RRS3} V. R\"odl, A. Ruci\'nski and E. Szemer\'edi, A note on perfect
matchings in uniform hypergraphs with large minimum collective
degree, \emph{Commentationes Mathematicae Universitatis Carolinae},
\textbf{49} (2008), 633--636.

\bibitem{RRS4} V. R\"odl, A. Ruci\'ski and E.
Szemer\'edi, Perfect matchings in large uniform hypergraphs with
large minimum collective degree, \emph{J. Comb. Theory Ser. A},
\textbf{116} (2009),  613--636.



\bibitem{RM} K. Markstr\"om and A.Ruci\'nski, Perfect matchings and Hamilton cycles
in hypergraphs with large degrees, \textit{Europ. J. Combin.}, \textbf{32} (2011), 677--687.


\bibitem{Tutte} W.T. Tutte, The factorization of linear graphs,\emph{ J. London Math. Soc.}, \textbf{22}
(1947),107--111.





\end{thebibliography}
\end{document}